\documentclass[11pt]{article} 
\usepackage{amssymb}

\newtheorem{thm}{Theorem}
\newtheorem{lem}{Lemma}

\newtheorem{cor}{Corollary}

\newtheorem{exm}{Example}

\begin {document}

\begin {center}
{\bf\Large On categorical semigroups}
\end {center}

\begin {center}
A.\,Kostin (Kiev, Ukraine)\\
B.\,Novikov (Kharkov, Ukraine)
\end {center}

\begin {abstract}
The structure of categorical at zero semigroups is studied from
the point of view their likeness to categories.
\end {abstract}

\bigskip

{\bf 1.} The connection between small categories and semigroups is
well-known. Let $\underline{C}$ be a small category. We join an
extra element $0$ to the set of its morphisms. The obtained set
$S(\underline{C})$ becomes a semigroup (with the zero $0$) with
respect to the operation $\ast:$
$$
f\ast g = \left\{ {\begin{array}{cl}
   fg & \mbox{if $fg$ is defined,}\\
   0 & \mbox{otherwise.}
    \end{array} } \right.
$$
Besides, $S(\underline{C})$ satisfies a condition:

\medskip

\centerline{if $f\ast g\ast h=0 $, then either $f\ast g=0 $ or
$g\ast h=0$.}

\medskip

Semigroups with this property are called {\it categorical at
zero}. We shall name them {\it $K$-semigroups} for brevity.

So, each small category is a $K$-semigroup. Of course, the
converse is not right. Moreover, there are natural examples of
$K$-semigroups, close to categories, but not being those. Here are
two of them:
\begin{exm}\label{exm1}{\rm
Any family of metric sets and all contracting mappings between
them.\footnote {We are indebted to Prof. S.\,Favorov for this
example.}
}\end {exm}
\begin{exm}\label{exm2}{\rm
 Any family of sets of the given infinite cardinal number $p$ and
all their mappings $f:A\to B$ such, that $|B\setminus fA|=p$.
}\end {exm}

Both examples do not contain identity morphisms and therefore are
not categories.

Here it is worth to mention also quasi-categories of C.\,Ehresmann
\cite{ehr}. One more example will be considered below.

This note is devoted to the study of the structure of
$K$-semigroups. We hope that it will be useful for understanding
what could be "near relatives" of categories.

For categories we use terminology from \cite {mcl}. Necessary data
from the theory of semigroups can be found in \cite{c-p} and
\cite{how}.

\medskip

{\bf 2.} We shall need two classes of semigroups.

Recall that a semigroup $S$ with a zero is called {\it
$n$-nilpotent} if $S^n=0$.

Now let  $G$ be a group, $I$ and $\Lambda$ be sets, $W=(w_{\lambda
i})_{\lambda \in \Lambda,\,i \in I}$ a matrix whose elements
$w_{\lambda i}$ are taken from the group with a jointed zero
$G\cup 0$. The set $(I\times G \times \Lambda)\cup 0$ with the
(associative) operation
$$
(i,g,\lambda)(j,h,\mu)=(i,gw_{\lambda j}h,\mu)
$$
is called {\it a Rees semigroup} and is denoted by
$M=M^0(G;\,I,\,\Lambda;\,W)$.

We shall use the special case of Rees semigroups, when $G=1$, and
write the non-zero elements of $M=M^0(1;\,I,\,\Lambda;\,W)$ as
$(i,\lambda)$. Then the multiplication has a form
$$
(i,\lambda)(j,\mu) = \left \{\begin{array}{lc}
   (i,\mu) &\mbox {\ if \ } w_{\lambda j}=1,\\
   0 &\mbox{\ if \ } w_{\lambda j} = 0.
 \end{array}\right.
$$

It is easy to see that except small categories, also 2-nilpotent
semigroups and Rees semigroups are
$K$-semigroups.\footnote{Generally speaking, it is possible to
turn every semigroup into a $K$-semigroup, if we join a zero
element to it, but such point of view, certainly, is ineffective.}

{\bf 3.} Further $S$ will denote a $K$-semigroup. We call the
subset ${\rm Ann}_l S = \{a \in S\mid aS = 0\}$ by the {\it left
annihilator} of $S$; similarly the {\it right annihilator} is the
subset ${\rm Ann}_r S = \{a\in S\mid Sa = 0\}$. The union
$$
{\rm Ann}_q S = {\rm Ann}_l S \bigcup {\rm Ann}_r S.
$$
is called a {\it quasi-annihilator}. Obviously, both left and
right annihilators (and, hence, the quasi-annihilator) are
two-sided ideals.

First we consider 3-nilpotent $K$-semigroups.

\begin{lem}\label{lem9:1}
$S$ is 3-nilpotent if and only if it coincides with ${\rm
Ann}_qS$.
\end {lem}

{\bf Proof.} Let $S$ is 3-nilpotent and $a \notin {\rm Ann}_lS$.
Then $ab\ne 0$ for some $b \in S$. Since $xab=0$ for all $x \in
S$, $xa=0$, i.\,e. $a \in {\rm Ann}_rS$. Hence $S = {\rm Ann}_qS$.

Conversely, if $S^3\ne 0$, then $abc\ne 0$ for some $a,b,c \in S$.
Then $ab \ne 0$ and $bc \ne 0$, i.\,e. $b \notin {\rm Ann}_r S\cup
{\rm Ann}_l S$, whence $S \ne {\rm Ann}_q S$. $\blacksquare$

\begin {cor} \label {cor9:0}
The quasi-annihilator of a $K$-semigroup is 3-nilpotent.
$\blacksquare$
\end {cor}

Lemma \ref{lem9:1} allows to build all 3-nilpotent $K$-semigroups.
Namely, let a set $A$ be given with a fixed element 0, two subsets
$B, C\subseteq A$ such, that $A=B\cup C$, \,$B\cap C \ni 0$, and a
mapping $\varphi: (B\setminus C) \times (C\setminus B) \to B \cap
C$, which satisfy conditions:

\begin {quote}
a) for every $b \in B\setminus C$ there is $c\in C\setminus B$
such that $\varphi (b, c) \ne 0$;
\end {quote}

\begin {quote}
b) for every $c \in C\setminus B$ there is $b \in B\setminus C$
such that $\varphi (b, c) \ne 0$.
\end {quote}

Define a multiplication on $A$:
$$
xy = \left \{{\begin {array} {lc}
   0, \ \mbox {if \ } x\in C\ \mbox {or \ } y\in B, \\
   \varphi (x, y), \ \mbox {if \ } x\in B\setminus C\
   \mbox {and \ } y\in C\setminus B. \\
 \end {array}} \right.
$$

If $xy \ne 0\ne yz$ then $y \in B\cap C$ and we get the
contradiction. Therefore for all $x,y,z\in A$ either $xy=0$ or
$yz=0$. It follows from here both 3-nilpotency and categoricity at
zero. Besides the conditions a) and b) provide equalities $C =
{\rm Ann}_l A, \ B = {\rm Ann}_r A$.

Now it is reasonable to consider the quotient semigroup $S/{\rm
Ann}_q S$. It turns out that $S$ is a ``splitting'' ideal
extension of ${\rm Ann}_q S$ in the following sense:

\begin {lem} \label {lem9:2}
The subset $T=(S \setminus {\rm Ann}_q S) \cup 0$ is a
subsemigroup.
\end {lem}

{\bf Proof.} Suppose that $a,b\in T$, $ab \in {\rm Ann}_q S$ and
$ab\ne 0$. Let, e.\,g., $ab \in {\rm Ann}_l S$. Then $abx=0$ for
all $x\in S$, whence $bx=0$, $b \in {\rm Ann}_l S$. $\blacksquare$

Note that $S /{\rm Ann}_q S\cong T$.

\medskip

{\bf 4.} To study the structure of $T$ we introduce the following
relations\footnote{They were defined by L.\,Gluskin \cite{glu} for
0-simple semigroups.} on $S$:

\medskip

 $\mathcal{P} = \{(a,b) \in S \times S\mid \forall x \in S\quad
xa = 0 \Leftrightarrow xb = 0\},$

 $\mathcal{Q} = \{(a,b) \in S \times S\mid \forall x \in S \quad
ax = 0 \Leftrightarrow bx = 0\},$

 $\mathcal{N} = \mathcal{P} \bigcap \mathcal{Q}$.

\medskip

In what follows we shall consider the restrictions of
$\mathcal{P,Q}$ and $\mathcal{N}$ on $T$ an denote them by the
same letters.

Obviously, $\mathcal{P,Q}$ and $\mathcal{N}$ are equivalences.
Furthermore, for each of them 0 forms the single-element class
(for example, if $a\in T$, $a\mathcal {P} 0$ then $Ta=0$, i.\,e.
$a \in {\rm Ann}_r S \cap T=0$).

Denote by $P_i$ and $Q_\lambda$  $\mathcal{P}$- and $\mathcal{Q}$-
classes respectively ($i \in I, \ \lambda \in \Lambda$, where $I$
and $\Lambda$ are sets of indexes). Set $N_{i\lambda}=P_i\cap
Q_\lambda$. So defined class $N_{i\lambda}$ is either empty or a
$\mathcal{N}$-class. Let
$$
P_i^0 =P_i\cup 0, \qquad Q_\lambda^0 =Q_\lambda \cup 0, \qquad
N_{i\lambda}^0 =N_{i\lambda} \cup 0.
$$

A homomorphism, for which the complete preimage of zero is
one-element, is called {\it 0-restricted}. A congruence,
corresponding to a 0-restricted homomorphism, we shall name also
{\it 0-restricted}.

\begin{lem}\label{lem9:3}
$\mathcal{N}$ is the greatest 0-restricted congruence on $T$.
\end{lem}

{\bf Proof.} Let $(a,b) \in \mathcal{N}, \ t \in T$. Show that
$(ta,tb) \in \mathcal{N}$. If $ta=0$ then also $tb=0$ (because
$(a,b) \in \mathcal{P}$), i.\,e. $(ta,tb) \in\mathcal{P}$. Let $ta
\ne 0 \ne tb$. If $xta=0$ then $xt=0$, whence $xtb=0$ and $(ta,tb)
\in \mathcal{P}$. The same way $(ta,tb) \in \mathcal{Q}$.

So $(ta,tb) \in \mathcal{N}$, similarly $(at,bt) \in \mathcal{N}$;
hence, $\mathcal{N}$ is a congruence. Evidently, $\mathcal{N}$ is
0-restricted.\footnote{Indeed $\mathcal{N}$ is also a congruence
on $S$, but not 0-restricted.}

Let $\rho$ be a 0-restricted congruence on $T$, $(a, b) \in \rho$.
If $xa=0$, then $(0,xb) \in \rho$, whence $xb=0$. Analogously,
$ax=0$ if and only if $bx=0$. Hence, $(a, b) \in \mathcal{N}$ and
$\mathcal{N}$ is the greatest 0-restricted congruence.
$\blacksquare$

\medskip

Now we elucidate, what is the quotient semigroup $T/\mathcal {N}$:

\begin {lem} \label {lem9:4}
$P_i Q_\lambda \subseteq N_{i\lambda}^0$.
\end {lem}

{\bf Proof.} Since $P_i^0$ is a right ideal and $Q_\lambda^0 $ is
a left one, then $P_i Q_\lambda \subseteq (P_i \cap Q_\lambda)
\cup 0 = N_{i\lambda}^0$. $\blacksquare$

\begin {cor} \label {cor9:1}
If $N_{i\lambda}$ and $N_{j\mu}$ are non-empty then $N_{i\lambda}
N_{j\mu} \subseteq N_{i\mu}^0$. $\blacksquare$
\end {cor}

\begin {lem} \label {lem9:5}
For every $i\in I, \ \lambda \in \Lambda$ either $0 \notin
Q_\lambda{}P_i$ or $Q_\lambda P_i=0$.
\end {lem}

{\bf Proof.} Suppose that $0 \in Q_\lambda P_i$, i.\,e. $yx=0$ for
some $x \in P_i$, $y \in Q_\lambda$. Choose arbitrary $u \in P_i,
\ v\in Q_\lambda$. From $(x,u) \in \mathcal{P}$ it follows $yu=0$;
from here and from $(y, v) \in \mathcal{Q}$ it follows $vu=0$.
Therefore $Q_\lambda P_i=0$. $\blacksquare$

\begin{cor}\label {cor9:2}
Assume that $N_{i \lambda}$ and $N_{j\mu}$ are non-empty. If
$Q_\lambda P_j=0$ then $N_{i\lambda} N_{j\mu} =0$. If $Q_{\lambda}
P_j \ne 0$ then $N_{i\lambda} N_{j\mu} \subseteq N_{i\mu}$.
$\blacksquare$
\end{cor}

In particular, we obtain some information about
$\mathcal{N}$-classes:

\begin{cor}\label {cor9:3}
$N_{i\lambda}^0$ is either a semigroup with zero multiplication or
a semigroup with an joined extra zero. $\blacksquare$
\end{cor}

Consider a Rees semigroup $M = M^0(1;\,I,\,\Lambda;\,W)$ with a
sandwich matrix $W=(w_{\lambda i})_{\lambda \in \Lambda,\,i \in
I}$, obeyed the condition:
$$
w_{\lambda i} = \left \{\begin{array}{cc}
   1 \mbox{\ if \ } Q_\lambda P_i \ne 0,\\
   0 \mbox{\ if \ } Q_\lambda P_i = 0.
 \end{array} \right.
 $$

Set $\varphi(0)=0$ and $\varphi (N_{i\lambda}) = (i,\lambda)$ for
every nonempty $\mathcal{N}$-class $N_{i\lambda}$. In that way a
mapping $\varphi:T/\mathcal{N} \to M$ is defined. By Corollary
\ref{cor9:2}
$$
\varphi (N_{i\lambda} N_{j\mu}) = \varphi (N_{i\mu}) = (i,\mu) =
(i,\lambda)(j,\mu) = \varphi (N_{i\lambda}) \varphi (N_{j\mu}),
$$
if $Q_\lambda P_i \ne 0 $, and
$$
\varphi (N_{i\lambda} N_{j\mu}) = \varphi (0) =0=(i,\lambda)(j,
\mu) = \varphi (N_{i\lambda}) \varphi (N_{j\mu}),
$$
if $Q_\lambda P_i = 0$. Hence, $\varphi$ is a monomorphism.

So we have proved the following assertion:

\begin{thm}\label{thm9:2}
Every $K$-semigroup $S$ is a splitting ideal extension of a
3-nil\-potent ideal $A$ by the subsemigroup $T=(S\setminus A)\cup
0$. The quotient semigroup of $T$ by its greatest 0-restricted
congruence is isomorphic to a subsemigroup of a Rees semigroup
with the one-element basic group. $\blacksquare$
\end{thm}

\medskip

{\bf 5.} The obtained results are interpreted for categories as
follows. Let $K$-semigroup $S$ is a small category. Then ${\rm
Ann}_qS=0$, i.\,e. $S=T$. Further, the fact that elements from
$\mathcal{P}$-class $P_i$ are annihilated at the left by the same
elements, means that $P_i$ is the set of all arrows which end in
the object $i$. Similarly, $Q_i$ is the set of arrows starting out
$i$. From here it follows that $I$ and $\Lambda$ can be identified
each with other and with the set of the objects of the category
$S$. It is easy to see that $Q_j P_i\ne 0 $ if and only if $i=j$.
Besides, $N_{ij} = {\rm Mor}(j,i)$, hence $N_{ij} N_{kl} = 0$ for
$j\ne k$, and $N_{ii}$ is a monoid. At last, the homomorphism
$\varphi$ is a functor from $S$ to $S/\mathcal{N}$, bijective on
the objects.

\medskip

We finish the article with an example of a ``non-category''.
\begin{exm}\label{exm3}{\rm
Let $\underline{C}$ be a small category, $\underline{D}$ and
$\underline{\Delta}$ its subcategories. For the simplicity we
assume that for every object $a$ from $\underline{D}$
$$
{\rm Mor}(a,\underline{\Delta})=\bigcup_{\alpha\in {\rm
Ob}\underline{\Delta}}{\rm Mor}(a,\alpha)\ne\emptyset,
$$
and choose a morphism $\varepsilon_a\in{\rm
Mor}(a,\underline{\Delta})$ for each $a\in{\rm Ob}\underline{D}$.
Let $\varepsilon_a:a\to\overline{a}$.

Denote
$$
{\rm Mor}(\underline{\Delta},\underline{D})= \bigcup_{\alpha\in
{\rm Ob}\underline{\Delta},\ a\in {\rm Ob}\underline{D}}{\rm
Mor}(\alpha,a)
$$
and $S={\rm Mor}(\underline{\Delta},\underline{D})\cup \{0\}$
where $0$ is an extra zero element. Define an operation $\ast$ on
$S$: for $f\in{\rm Mor}(\alpha,a)$, $g\in{\rm Mor}(\beta,b)$
$$
g\ast f = \left\{ {\begin{array}{ll}
   g\varepsilon_af, & \mbox{if $\overline{a}=\beta$,}\\
   0 & \mbox{otherwise.}
    \end{array} } \right.
$$
It is easy to see that $S$ becomes a $K$-semigroup. In this case
${\rm Ann}_r S=0$ and ${\rm Ann}_l S=\{g\mid {\rm dom}g\in{\rm
Ob}\underline{\Delta}\setminus\varepsilon(\underline{D})\}$ where
$\varepsilon(\underline{D})=\{\overline{a}\mid a\in {\rm
Ob}\underline{D}\}$.
}\end{exm}


\end {document}